\documentclass[mathpazo]{cmr}

\renewcommand\thisnumber{x}
\renewcommand\thisyear {202x}
\renewcommand\thismonth{x}
\renewcommand\thisvolume{x}

\usepackage{amsmath,amsfonts,amssymb, amsthm, xfrac}
\usepackage{enumerate}
\usepackage{mathrsfs}
\usepackage{graphicx}
\usepackage{tikz}
\usepackage{amsmath}
\usepackage{amscd}

\newcommand{\N}{\mathcal{N}}
\newcommand{\Z}{\mathbb{Z}}

\newcommand{\C}{\mathcal{C}}

\newcommand{\Fq}{\mathbb{F}_q}

\newcommand{\ord}{\mathop{\rm ord}}

\newcommand{\lcm}{\mathop{\rm lcm}}
\newcommand{\cyc}{\mathop{\rm Cyc}}

\newcommand{\G}{\mathscr{G}}

\newcommand{\Ex}{\mathbb{E}}
\newcommand{\U}{\mathcal{U}}
\newcommand{\Ap}{\text{AP}}
\renewcommand{\L}{\mathcal{L}}
\newcommand{\Gc}{\text{GCP}}
\renewcommand{\O}{\mathcal{O}}

\title{How random are affine maps over finite cyclic groups and generalized cyclotomic permutations over finite fields?}

\author{Hugo Teixeira} 
\address{Johann Radon Institute for Computational and Applied Mathematics, Austrian Academy of Sciences,
Altenbergestr. 69, 4040 Linz, Austria}
\email{hugo.teixeira@ricam.oeaw.ac.at}



\begin{document}

\begin{abstract}
  We analyze the dynamical behavior of affine permutations over the ring $\Z_{p^e}$ and generalized cyclotomic permutations (GCPs) over the finite field $\Fq$. We use the description of their functional graphs to determine two statistical parameters: the expected number of cycles and the expected length of a cycle that contains a random element. We obtain explicit closed-form expressions for both expectations and an exact reduction from GCPs to affine permutations. We compare these parameters with the corresponding expectations for random permutations, to determine how closely the cycle statistics of each family resemble those of random permutations. 
\end{abstract}

\maketitle
\section{Introduction}

A \textit{discrete dynamical system} is a pair $(S,f)$, where the state space $S$ is a finite set and $f$ is a self map of $S$. The iteration of the function $f$ determines the evolution of each state of the state space $S$ in a discrete time \cite{martelli1999}. Associated with the discrete dynamical system $(S,f)$ is a directed graph  $\G_f(S)$, known as the \textit{functional graph} of $f$ over $S$, with vertex set $V=S$ and directed edges $E=\{(x,f(x))\mid x\in S\}$. This digraph is a visual representation of the evolution of each state of $S$ under the function $f$.

Flajolet and Odlyzko \cite{Fla} determined the asymptotic behavior of the functional graph associated with the discrete dynamical system $(S,f)$, where $f$ is a function selected uniformly at random. They obtain stochastic parameters, such as the expectation for cycle length, number of connected components and number of cycles of a given length for a function taken with uniform probability from the set $S^S$ of all self maps of $S$. The same statistical properties for random permutations can be found in \cite{Stanley_2011}. Therefore, obtaining such parameters for a smaller set of functions is a natural way to compare their dynamical behavior to those of a truly random function. Many function families over finite structures have been investigated in that sense \cite{PanRei,QurPan1,QurPan2,MaPaQu}.

The case where the state space $S$ is a finite field is particularly important due to its diverse applications in areas such as cybersecurity \cite{NySo,BCP}, pseudorandom number generation \cite{Ei-Eh,NiWi15,ToWi07,Winterhof10} and biology \cite{BBK, KWV}. Within this setting, dynamical systems induced by permutations are often sought out for being building blocks in cryptographic protocols. In this case, the functional graph consists exclusively of disjoint cycles, so that its dynamical behavior is completely determined by its cycle structure. Consequently, the number and lengths of cycles provide natural parameters for studying and comparing families of permutations over finite fields.

In this paper, we investigate stochastic properties of the functional graphs of two families of permutations: affine permutations over the finite cyclic group of order $p^e$ and \textit{generalized cyclotomic permutations} over finite fields. Introduced by Niederreiter and Winterhof \cite{NW05} and Wang \cite{Wan07}, \textit{generalized cyclotomic mappings} (or GCMs) of index $d$ are a generalization of monomial maps. Let $q$ be a prime power and $d$ a divisor of $q-1$. An index $d$ GCM of the finite field $\Fq$ is any function $f:\Fq\to\Fq$ that satisfies $f(0)=0$ and agrees with a monomial function when restricted to each of the $d$ cosets of the unique index $d$ subgroup of $\Fq^*$. A generalized cyclotomic permutation (or GCP) is a bijective GCM. The functional graphs of index $d$ GCMs are described by Bors, Panario and Wang in \cite{ADS}.
 
Clearly, any monomial permutation over $\Fq$ is an index $1$ GCP. Their functional graph, described in \cite{deng2013isomorphic}, is highly dependent on the divisors of $q-1$. However, any permutation of $\Fq$ that fixes $0$ is a GCP of index $q-1$. Thus, a randomly selected GCP is essentially a random permutation of $\Fq^*$. A natural question is how the cycle statistics of GCPs change as the index increases, and how far from $q-1$ the index can be while maintaining cycle statistics comparable to those of uniformly random permutations.
 
 We obtain a group isomorphism through discrete logarithm between the index $d$ subgroup of $\Fq^*$ and the additive group $\Z_{s}$, where $s=\frac{q-1}{d}$. We then obtain the expected number of cycles and the expected length of the cycle that contains a random element in the functional graph of a randomly selected generalized cyclotomic permutation in terms of the same expectations for affine permutations over $\Z_s$. 
 
 Moreover, we determine explicit closed formulas for these expectations in affine permutations over $\Z_{p^e}$, where $p$ is a prime and $e$ is a positive integer. When $p^e$ divides $q-1$, we consider index $d=\frac{q-1}{p^e}$ GCPs. The exponent $e$ provides a natural parameter for controlling the size of $d$, and consequently its distance from the maximal possible index $q-1$.  We derive the expectations for GCPs of this particular index and observe their asymptotics as $e$ grows. Finally, we provide a comparison between these expectations with the corresponding expectations for uniformly random permutations over a set of the same size.
 
 
 We end this section by providing the structure of the paper. In Section \ref{GCM}, we describe the dynamical behavior of GCMs, determined by Bors, Panario and Wang in \cite{ADS}. In Section \ref{SGCP}, we derive a formula for the two statistical properties of index $d$ GCPs in terms of the corresponding properties for affine permutations over $\Z_s$. In Section \ref{Aff}, we present Deng's \cite{deng2013isomorphic} characterization of functional graphs of affine maps over the finite cyclic group $\Z_{s}$. Using this characterization, in Section \ref{SAff} we derive the stochastic parameters for affine permutations over $\Z_{p^e}$. In Section \ref{Random}, we provide a comparison between the properties from affine permutations over $\Z_{p^e}$ and GCPs against the corresponding expectations for random permutations. In Section \ref{Con}, we conclude the paper by summarizing our findings and presenting problems for future research.

\section{Preliminary properties}\label{Pre}
In this section we present a list of lemmas needed in the proofs of our main results. 

Throughout the paper, $p$ is a prime, $[a,b]$ is the set of nonnegative integers greater or equal to $a$ and smaller or equal to $b$, $\Z_n$ is the finite additive group of integers modulo $n$ and  $\ord_n(a)$ is the smallest positive integer $k$ such that $a^k\equiv1\pmod n$, for $a$ coprime with $n$. Moreover, we define $\nu_p(a)$ as the greatest $k\in\Z$ such that $p^k$ divides $a$ and set $\nu_p(0)=\infty$. 

\begin{lemma}[Lifting The Exponent Lemma \cite{Manea}] \label{lif-exp} 
	Let $x$ and $n$ be integers, and $p$ be a prime. 
    \begin{enumerate}
		\item If $p$ is an odd prime and $\nu_p(x-1)>0$, then $\nu_p({x^n-1})-\nu_p({x-1})=\nu_p(n)$.
		\item If $p=2$ and $n$ is odd, then $\nu_p(x^n-1)-\nu_p(x-1)=\nu_p(n)=0$.
		\item If $p=2$, $n$ is even and $x\equiv 1\pmod 4$, then $\nu_p(x^n-1)-\nu_p(x-1)=\nu_p(n)$.
		\item If $p=2$, $n$ is even and $x\equiv 3\pmod 4$, then $\nu_p(x^n-1)-\nu_p(x-1)=\nu_p(\frac{x+1}{2})+\nu_p(n)$.
	\end{enumerate}     
\end{lemma}

\begin{lemma}\label{ord}
    Let $p$ be a prime, $j$ a nonnegative integer, $k$ a positive integer and $i\in[1,p-1]$. If $p$ is odd, $d=\ord_p(i)$ and $t=\nu_p((jp+i)^d-1)$, then $\ord_{p^k}(jp+i)=d\cdot p^{\max\{0,k-t\}}$. 
    
    If $p=2$, then $ord_{2^k}(4j+1)=2^{\max\{0,k-2-\nu_2(j)\}}$. Moreover, $\ord_{2}(4j+3)=1$ and, for all $k>1$, $\ord_{2^k}(4j+3)=2^{\max\{1,k-2-\nu_2(j+1)\}}$. 
\end{lemma}
\begin{proof}
    Let $p$ be an odd prime. For $k=1$, we have that $jp+i\equiv i\pmod p$, thus $\ord_p(jp+i)=\ord_p(i)=d$. Moreover, since $p$ divides $p^k$ we have that $\ord_{p^k}(jp+i)=d\cdot r$, for some positive integer $r$. Hence, by Lemma \ref{lif-exp}, if $k>1$ we have that $$\nu_p((jp+i)^{\ord_{p^k}(jp+i)}-1)=\nu_p((jp+i)^{d\cdot r}-1)=\nu_p((jp+i)^{d}-1)+\nu_p(r)\geq k.$$ This implies that $r$ is the minimal positive integer that satisfies $\nu_p(r)\geq k-t$. Hence $r=p^{\max\{0,k-t\}}$.

   For $p=2$, from Lemma \ref{lif-exp}, we have that $ord_{2^k}(4j+1)$ is the minimal nonnegative integer $r$ that satisfies $$\nu_2((4j+1)^{r}-1)=\nu_2(4j)+\nu_2(r)\geq k.$$ Hence, we get that $ord_{2^k}(4j+1)=2^{\max\{0,k-2-\nu_2(j)\}}$. 
   
   For $4j+3$, we have trivially that $\ord_{2}(4j+3)=1$. If $k>1$, the multiplicative order of any odd integer in $\Z_{2^k}$ is a power of $2$, since it divides $2^{k-1}$, the cardinality of the multiplicative group of $\Z_{2^k}$. Moreover, we observe that $\nu_2((4j+3)-1)=1$, hence $\ord_{2^k}(4j+3)>1$. Therefore, from Lemma \ref{lif-exp} we get that $\ord_{2^k}(4j+3)$ is the minimal positive integer $r$ that satisfies $$\nu_2((4j+3)^{2^r}-1)=2+\nu_2(j+1)+r\geq k.$$ Thus, we conclude that $\ord_{2^k}(4j+3)=2^{\max\{1,k-2-\nu_2(j+1)\}}$.
\end{proof}
The following lemma is known as Hensel's lifting lemma. 
\begin{lemma}{\cite[Chapter 5, Theorem 3]{Koblitz}}\label{hensel}
    If $f(x)\in\Z[x]$ and $a\in\Z_p$ satisfies $f(a)\equiv0\pmod p$ and $f'(a)\not\equiv0\pmod p$, then there is a unique $\alpha_k\in\Z_{p^k}$ such that $f(\alpha_k)\equiv0\pmod{p^k}$ and $\alpha_k\equiv a\pmod p$.
\end{lemma}
\begin{lemma}\label{sumnup}
    If $p$ is a prime and $e>1$ a positive integer, then $$\sum_{j=1}^{p^{e}-1}\nu_p(j)=\frac{p^{e}-1}{p-1}-e.$$ 
\end{lemma}
\begin{proof}
    For any prime $p$ and $r\in[0,e-1]$, there are $p^{e-r}-p^{e-r-1}$ integers $j$ in $\Z_{p^{e-1}}$ with $\nu_p(j)=r$. Thus, we obtain $$\sum_{j=1}^{p^{e}-1}\nu_p(j)=\sum_{r=0}^{e-1}r(p^{e-r}-p^{e-r-1})=\frac{p^{e}-1}{p-1}-e.$$ 
\end{proof}

\section{Generalized Cyclotomic Mappings}\label{GCM}
We begin by summarizing some of the results presented in \cite{ADS} that determine the dynamical behavior of generalized cyclotomic mappings. Let $q$ be a prime power, $d$ a divisor of $q-1$ and $s=\frac{q-1}{d}$. For any primitive element $\xi$ of $\Fq^*$, let $C_0=\left\{\xi^{j\cdot d}\mid j\in \left[0,s-1\right]\right\}$ be the unique index $d$ subgroup of $\Fq^*$. For $i\in [0,d-1] $, let $C_i=\xi^iC_0$ be the $d$ cosets of $C_0$. A \textit{generalized cyclotomic mapping of $\Fq$ of index $d$} is any function $f$ that satisfies $$f(x)=\begin{cases}
    m_0x^{a_0}&\text{if $x\in C_0$},\\
    m_1x^{a_1}&\text{if $x\in C_1$},\\
    \hspace{.4cm}\vdots& \hspace{.5cm}\vdots\\
    m_{d-1}x^{a_{d-1}}&\text{if $x\in C_{d-1},$}\\
    0&\text{if $x=0$},
\end{cases}$$ where $m_i\in\Fq^*$ and $a_i\in[0,q-2]$. We have that $f$ induces a unique coset map $\overline{f}:[0,d-1]\to[0,d-1]$ on the indices of the $d$ cosets such that $f(C_i)\subseteq C_{\overline{f}(i)}$. We are interested in the case when $f$ is a generalized cyclotomic permutation (GCP). This implies that $m_i$ is nonzero for every $i\in [0,d-1]$ and the connected component generated that contains zero is always a cycle of length one (also known as a \textit{fixed point}). Thus, we omit the case $x=0$ and consider only the restriction of $f$ to $\Fq^*$ hereafter. 

Let $\Gc(d,\Fq^*)$ be the set of all generalized cyclotomic permutations of index $d$ over $\Fq^*$. If $f\in\Gc(d,\Fq^*)$, the induced coset map $\overline{f}$ of $f$ is a permutation of the set $[0,d-1]$. Therefore, hereafter we denote by $\sigma_f\in S_d$ the representation of $\overline{f}$ as the element of the symmetric group $S_d$ written as a product of disjoint cycles.

Writing $m_i=\xi^{b_i}$, for some integer $b_i \in [0,q-2]$, we rewrite $f$ as a function acting on the exponent of the element's representation as a power of the primitive element $\xi$. 
$$f(\xi^t)=\begin{cases}
    \xi^{a_0t+b_0}&\text{if $\xi^t\in C_0$},\\
    \xi^{a_1t+b_1}&\text{if $\xi^t\in C_1$},\\
    \hspace{.4cm}\vdots& \hspace{.5cm}\vdots\\
    \xi^{a_{d-1}t+b_{d-1}}&\text{if $\xi^t\in C_{d-1}$}.
\end{cases}$$ For each $i\in[0,d-1]$, we define $f_i(t)=a_it+b_i$.

\begin{proposition}\label{sigmaGCP}
Let $\sigma$ be any element of $S_d$. There are $(s\cdot\varphi(s))^d$ generalized cyclotomic permutations of index $d$ over $\Fq$ with induced coset map $\sigma$. 
\end{proposition}
\begin{proof}
Let $\sigma\in S_d$ and $f$ be a generalized cyclotomic permutation of index $d$ over $\Fq$ with induced coset map $\sigma_f=\sigma$. We observe that $f$ is a permutation only if it is a bijection between $C_i$ and $C_{\sigma(i)}$, for every $i\in[0,d-1]$. Each $f_i$ induces the map $$t\mapsto a_it+\frac{a_ii+b_i-\sigma(i)}{d}=a_it+c_i$$ in $\Z_s$, which is completely determined by $a_i$ and $c_i$ modulo $s$ and is bijective if and only if $\gcd(a_i,s)=1$. Therefore, we conclude that $f$ is a permutation if and only if each $a_i$ is invertible modulo $s$, and there are $(s\cdot\varphi(s))^d$ such maps.
\end{proof}
\begin{corollary}
    There are $d!\cdot (s\cdot \varphi(s))^d$  generalized cyclotomic permutations of index $d$.
\end{corollary}

\section{Stochastic properties of GCPs}\label{SGCP}
For any finite set $S$ and map $f$ defined over $S$, we denote by $\N(f,S)$ the number of cycles in $\G_f(S)$, $\L(f,S)$ the set of lengths of cycles in $\G_f(S)$ and, for $\ell\in\L(f,S)$, $\N_\ell(f,S)$ is the number of cycles of length $\ell$ in $\G_f(S)$. Finally, for $x\in S$ let $\ell_x(f,S)$ be the length of the cycle that contains $x$ in $\G_f(S)$. 

Let $\Ap(\Z_{s})$ be the set of all affine permutations over $\Z_s$. We now establish the precise relationship between cycles generated by elements in $\Ap(\Z_{s})$ and those generated by elements of $\Gc(d,\Fq^*)$. We determine the expected cycle length of a random element and expected number of cycles of a randomly selected index $d$ generalized cyclotomic permutation in terms those same properties for affine functions over $\Z_s$.

Let $f\in\Gc(d,\Fq^*)$ with coset map $\sigma_f$. We define $\C(\sigma_f)$ as the set of all cyclically ordered subsets $\Lambda=(\lambda_0,\ldots,\lambda_{n-1})$ of the integers $[0,d-1]$, for all $n\in[1,d]$, such that $(\lambda_0\ldots\lambda_{n-1})$ is a cycle in $\sigma_f$. For $\Lambda\in\C(\sigma_f)$, let $C_\Lambda=\bigcup_{\lambda\in\Lambda}C_\lambda$. 
If $n$ is a positive integer, we define the $n$th iteration of $f$ as $f^{(1)}(x)=f(x)$ and $f^{(n)}(x)=f(f^{(n-1)}(x))$, for $n>1$. 

In general, if $\G_g(S)$ and $\G_h(R)$ are functional graphs associated with the functions $g$ over $S$ and $h$ over $R$ respectively, then $\G_g(S)$ and $\G_h(R)$ are \textit{isomorphic graphs} if there is a bijection $\phi$ between $S$ and $R$ such that $\phi^{-1}(g(\phi(a)))=h(a)$, for any $a\in R$. This means that the graphs $\G_g(S)$ and $\G_h(R)$ are equal up to relabeling of vertices.


\begin{proposition}\label{sum}
    Let $f\in\Gc(d,\Fq^*)$ with coset map $\sigma_f$. Then
    $$\N(f,\Fq^*)=\sum_{\Lambda\in\C(\sigma_f)}\N(A_\Lambda,\Z_{s})$$ and $$\sum_{\ell\in\L(f,\Fq^*)}\ell^2\N_\ell(f,\Fq^*)=\sum_{\Lambda\in\C(\sigma_f)}|\Lambda|^2\sum_{\ell\in\L(A_\Lambda,\Z_s)} \ell^2\cdot\N_\ell(A_\Lambda,\Z_{s})$$ where $A_\Lambda(t)=\frac{f_{\lambda_{n-1}}\circ\cdots\circ f_{\lambda_0}(\lambda_0+td)-\lambda_0}{d}$ for any $\Lambda=(\lambda_0,\ldots,\lambda_{n-1})\in\C(\sigma_f)$.
\end{proposition}
\begin{proof}
Every cycle in the functional graph of $f$ is contained in a set $C_\Lambda$, for some $\Lambda\in\C(\sigma_f)$. If $\Lambda=(\lambda_0,\ldots,\lambda_{n-1})\in\C(\sigma_f)$, then the $n$th iteration of $f$ is a permutation of $C_{\lambda_0}$. Moreover, the bijection $\iota_{\lambda_i}:\Z_{s}\to C_{\lambda_i}$, $t\mapsto\xi^{\lambda_i+dt}$ creates a commutative diagram 
$$
\begin{CD}
C_{\lambda_0} @>f^{(n)}>> C_{\lambda_0} \\
@V\iota_{\lambda_0}^{-1}VV @AA\iota_{\lambda_0}A \\
 \Z_{s} @>A_\Lambda>>  \Z_{s}
\end{CD}
$$
where $A_\Lambda(t)=\frac{f_{\lambda_{n-1}}\circ\cdots\circ f_{\lambda_0}(\lambda_0+td)-\lambda_0}{d}.$ Hence, there is a one to one correspondence between the cycles generated by $A_\Lambda$ and the cycles generated by $f$ in $C_\Lambda$. In addition,  each cycle of length $\ell$ in the functional graph of $A_\Lambda$ corresponds to a cycle of length $n\cdot\ell$ in the functional graph of $f$. 

Moreover, for any circular permutation $(\lambda_j\ldots\lambda_{j-1})$ of $({\lambda_0\ldots{\lambda_{n-1}})}$, where the indices are taken modulo $n$, taking the bijection $\phi(t)=f_{\lambda_{j}}\circ\cdots\circ f_{\lambda_{n-1}}(\lambda_j+~td)$, we have that the functional graph of $$\frac{f_{\lambda_j}\circ\cdots\circ f_{\lambda_{j-1}}(\lambda_{j-1}+td)-\lambda_{j-1}}{d}$$ over $\Z_{p^e}$ is isomorphic to the graph of $\G_{A_\Lambda}(\Z_{s})$. Since this goes for every $\Lambda\in\C(\sigma_f)$, the result follows.  
\end{proof}

\begin{proposition}\label{sumgcpcycle}
    Let $d$ be a divisor of $q-1$ and set $s=\frac{q-1}d$. Let $g(t)=at+b\in\Ap(Z_{s})$, $n\in[1,d]$ and $\Lambda=(\lambda_0,\ldots,\lambda_{n-1})$ be a cyclically ordered subset of $[0,d-1]$. The number of generalized cyclotomic permutations of index $d$ that generate a cycle in the cosets with index in $\Lambda$ and have induced map $A_\Lambda(t)=g(t)$ is $$(d-n)!(s\varphi(s))^{d-1}.$$
\end{proposition}
\begin{proof}
    For fixed $\Lambda=(\lambda_0\ldots\lambda_{n-1})$ and $A_\Lambda$, the map $A_\Lambda$ is determined by the $n$ maps $f_{\lambda_i}(t)$, for $\lambda_i\in\Lambda$. Choosing $n-1$ of them arbitrarily among the $s\varphi(s)$ choices, the last one is uniquely determined modulo $s$ by the fixed map $A_\Lambda$. For the remaining $d-n$ cosets, there are $(d-n)!$ choices for the coset permutation and $s\varphi(s)$ choices for each of the $d-n$ maps.
\end{proof}

\begin{theorem}\label{expecGCP}
    Let $q$ be a prime power, $d$ a divisor of $q-1$ and set $s=\frac{q-1}{d}$. If $f$ and $x$ are selected uniformly at random from the sets $\Gc(d,\Fq^*)$ and $\Fq^*$, respectively, then 
     $$\Ex[\N(f,\Fq^*)]=H_d\cdot\Ex[\N(g,\Z_{s})]\quad\text{and}\quad\Ex[\ell_x(f,\Fq^*)]=\frac{d+1}{2}\Ex[\ell_t(g,\Z_{s})]$$ where $g$ and $t$ are uniformly distributed over $\Ap(\Z_{s})$ and $\Z_s$, respectively, and $H_d$ is the $d$th harmonic number.
\end{theorem}
\begin{proof}
    There are $\displaystyle\binom{d}{n}(n-1)!$ cyclically ordered subsets of $[0,d-1]$, for each $n\in[1,d]$. From Propositions \ref{sum} and \ref{sumgcpcycle}, each cyclically ordered subset $\Lambda$ and fixed affine function $A_\Lambda$ occurs on $(d-n)!\cdot (s\cdot\varphi(s))^{d-1}$ GCPs. Thus, we have that\begin{align*}
        \sum_{f\in\Gc(d,\Fq^*)}\N(f,\Fq^*)&=\sum_{n=1}^d\binom{d}{n}(n-1)!(d-n)!(s\varphi(s))^{d-1}\sum_{g\in\Ap(\Z_{s})}\N(g,\Z_s)\\
        &=d!(s\varphi(s))^{d-1}H_d\sum_{g\in\Ap(\Z_{s})}\N(g,\Z_s).
    \end{align*}Dividing by $|\Gc(d,\Fq^*)|=d!(s\varphi(s))^d$, the first expectation follows.

    Moreover, we have that \begin{align*}    \sum_{f\in\Gc(d,\Fq^*)}\sum_{\ell\in\L(f,\Fq^*)}\ell^2\N_\ell(f,\Fq^*)\\
    &\hspace{-4cm}= \sum_{f\in\Gc(d,\Fq^*)}\sum_{\Lambda\in\C(\sigma_f)}|\Lambda|^2\sum_{\ell\in\L(A_\Lambda,\Z_s)} \ell^2\cdot\N_\ell(A_\Lambda,\Z_{s})\\
    &\hspace{-4cm}=\sum_{n=1}^d\binom{d}{n}(n-1)!(d-n)!(s\varphi(s))^{d-1}n^2\sum_{g\in\Ap(\Z_s)}\sum_{\ell\in\L(g,\Z_s)} \ell^2\cdot\N_\ell(g,\Z_{s})\\
    &\hspace{-4cm}=d!(s\varphi(s))^{d-1}\sum_{n=1}^dn\sum_{g\in\Ap(\Z_s)}\sum_{\ell\in\L(g,\Z_s)} \ell^2\cdot\N_\ell(g,\Z_{s})\\
    &\hspace{-4cm}=\frac{d(d+1)}{2}d!(s\varphi(s))^{d-1}\sum_{g\in\Ap(\Z_s)}\sum_{\ell\in\L(g,\Z_s)} \ell^2\cdot\N_\ell(g,\Z_{s}).
    \end{align*}    Dividing by $|\Fq^*|\cdot|\Gc(d,\Fq^*)|=(q-1)d!(s\varphi(s))^{d}$ the second expectation follows.
\end{proof} 

This determines the precise relation between the expectations for affine maps over $\Z_s$ and for $GCP$s over $\Fq^*$. In the next sections, we derive a closed-form expression for these expectations on affine maps over $\Z_s$, when $s$ is a prime power and provide a comparison with random permutations.

\section{Affine functions over finite cyclic groups}\label{Aff}

The dynamical behavior of any nonzero monomial function over $\Fq^*$ is isomorphic via discrete logarithm to the behavior of an affine function over the additive group $\Z_{q-1}$. More precisely, if $\xi$ is a primitive element of $\Fq^*$, $a\in[0,q-1]$ and $m=\xi^b\in\Fq^*$, then the functional graphs of $f(x)=mx^a$ over $\Fq^*$ and $g(x)=at+b$ over $\Z_{q-1}$ are isomorphic through the bijection $\phi:\Z_{q-1}\to\Fq^*$ given by $\phi(t)=\xi^t$.

The functional graph of affine functions over finite cyclic groups is completely described by Deng in \cite{deng2013isomorphic}. In this section, we summarize some of his  findings.

\begin{definition}{\cite[Formula (1)]{deng2013isomorphic}}\label{decomposition}
Let $\G_1,\ldots,\G_k$ be directed graphs. We define $\prod_{i=1}^k\G_i$ as the directed graph with vertex set $V=\{(v_1,\ldots,v_k)\mid v_i\in\G_i, i\in[1,k]\}$ and edges from $(v_1,\ldots,v_k)$ to $(w_1,\ldots,w_k)$ if and only if $v_i$ is directed to $w_i$ in $\G_i$, for all $i\in[1,k]$.      
\end{definition}

For a family of maps $g_i:S_i\to S_i$, with $i\in[1,k]$, we define $\bigotimes_{i=1}^k g_i:\prod_{i=1}^k S_i\to \prod_{i=1}^k S_i$ as the map $(a_1,\ldots,a_k)\mapsto(g_1(a_1),\ldots,g_k(a_k))$. Then, from the previous definition, its functional graph is given by  $\G_{\bigotimes_{i=1}^k g_i}(\prod_{i=1}^k S_i)=\prod_{i=1}^k\G_{g_i}(S_i)$. 

\begin{definition}
    For any positive integer $n$, we denote by $\cyc(n)$ the directed graph composed of a cycle with $n$ vertices. The directed graph $\G=\bigoplus_{i=0}^kr_i\cdot\cyc(n_i)$ is the graph composed by $r_i$ copies of $\cyc(n_i)$, for any nonnegative integers $r_i$ and positive integers $n_i$, $i\in[0,k]$.
\end{definition}
\begin{lemma}{\cite[Lemma 1]{deng2013isomorphic}}\label{lcm}
    If $n_1$ and $n_2$ are positive integers,  the graph $\cyc(n_1)\times\cyc(n_2)$ is isomorphic to $\gcd(n_1,n_2)\cdot\cyc(\lcm(n_1,n_2))$.
\end{lemma}

Moreover, for affine maps over $\Z_s$ we have the following.

\begin{lemma}{\cite[Theorem 2]{deng2013isomorphic}}\label{Theorem2}
    Let $s=p_1^{e_1}\cdots p_k^{e_k}$ and $f(x)=ax+b$, for some $a,b\in\Z$. Then $\G_f(\Z_s)=\prod_{i=1}^k\G_{f}(\Z_{p_i^{e_i}})$. 
\end{lemma}

This lemma combines Definition \ref{decomposition} with the Chinese Remainder Theorem and allows us to determine the behavior of $\G_f(\Z_s)$ by studying the behavior of $\G_{f}(\Z_{p_i^{e_i}})$, for every $p_i^{e_i}$ in the prime decomposition of $s$. 

We recall that when $f$ is a permutation, its functional graph is composed of cycles with no attached trees. Hence, its description is given by the number of cycles of each length.

\begin{theorem}{\cite[Lemmas 3 and 7]{deng2013isomorphic}}\label{affinecycles}
    Let $p$ be a prime and $e$ be a positive integer. For $a,b\in\Z$, the affine function $f(x)=ax+b$ is a permutation over $\Z_{p^e}$ if and only if $\gcd(a,p)=1$. In that case, the functional graph of $f$ over $\Z_{p^e}$ is given by the following two cases.
    \begin{enumerate}
        \item If $\nu_p(a-1)\leq \nu_p(b)$, then $$\G_f(Z_{p^e})=\bigoplus_{i=1}^{\ord_{p^e}(a)}n_i\times \cyc(i),$$ where $n_i=\frac{1}{i}\sum_{d\mid i}\mu(d)\gcd(p^e,a^{\frac{i}{d}}-1)$.
        \item If $\nu_p(a-1)>\nu_p(b)$, then $$\G_f(Z_{p^e})=p^{e-s}\times\cyc(p^s),$$ where $$p^s=
        \begin{cases}
            p^{e-\nu_p(b)}&\text{if $p\neq 2$ or $\nu_p(a-1)>1$}\\
            2 &\text{if $p= 2$, $\nu_p(a-1)=1$ and $e+1\leq\nu_2(a^2-1)$}\\
            2^{e+2-\nu_p(a^2-1)}&\text{if $p= 2$, $\nu_p(a-1)=1$ and $e+1>\nu_2(a^2-1)$}
        \end{cases}.$$
    \end{enumerate}
\end{theorem}

In the first case, we have that $n_i$ is zero for some choices of $i\in[1,\ord_{p^e-1}(a)]$. The following result by Sha provides an alternative way to determine the cycle lengths in the functional graph of $f$ when $\nu_p(a-1)\leq \nu_p(b)$.

\begin{theorem}{\cite[Proposition 2.6]{Sha2012Digraphs}}\label{altaffinecycles}
    Let $p$ be a prime and $e$ be a positive integer. Let $f(x)=ax+b$ be such that $\gcd(a,p)=1$ and $\nu_p(a-1)\leq \nu_p(b)$. Then 
    $$\G_f(\Z_{p^e})=\bigoplus_{k=0}^{e}\frac{\varphi(p^k)}{\ord_{p^k}(a)}\times \cyc(\ord\nolimits_{p^k}(a))$$
 and the total number of cycles in $\G_f(\Z_{p^e})$ is $\sum_{j=0}^e\frac{\varphi(p^j)}{\ord_{p^j}(a)}$, where $\ord_{p^0}(a)=1$.
\end{theorem}
The drawback of this notation is that we might have $\ord_{p^k}(a)=\ord_{p^r}(a)$ for some $k\neq r$.

\section{Stochastic properties of affine maps over $\Z_{p^e}$}\label{SAff}

Although there is extensive research on the dynamical behavior of the maps in $\Ap(\Z_{s})$, including the complete characterization of cycle and tree decomposition, in this section we provide new closed-form expressions for two stochastic properties of their functional graphs in the case where $s$ is a prime power. 

Let $p$ be a prime and $e$ a positive integer. The expectation for the number of cycles of $f\in\Ap(\Z_{p^e})$ selected uniformly at random is determined by $$\Ex[\N(f,\Z_{p^e})]=\frac{1}{|\Ap(\Z_{p^e})|}\sum_{f\in\Ap(\Z_{p^e})}\N(f,\Z_{p^e}).$$ 
The expectation for the length of the cycle that contains $x\in\Z_{p^e}$ in the functional graph of $f\in\Ap(\Z_{p^e})$, where $x$ and $f$ are selected uniformly at random, is determined by $$\Ex[\ell_x(f,\Z_{p^e})]=\frac{1}{p^e\cdot|\Ap(\Z_{p^e})|}\sum_{f\in\Ap(\Z_{p^e})}\sum_{\ell\in\L(f,\Z_{p^e})}\ell^2\N_\ell(f,\Z_{p^e}).$$ Since the number of affine permutations over $\Z_{p^e}$ is $p^e\cdot\varphi(p^e)=p^{2e-1}(p-1)$, to calculate the expectations we need to determine the value for each of the sums. We begin with the case where $p=2$.

\subsection{Stochastic properties for $p=2$}

\begin{theorem}\label{cyc2e}
     $$\sum_{f\in\Ap(\Z_{2^e})}\N(f,\Z_{2^e})=(e^2+3e+8)2^{2e-4}.$$ 
\end{theorem}
\begin{proof}
    We begin with the trivial case $e=1$, where there are only two affine maps: the identity and $x\mapsto x+1$. The first generates two cycles and the second generates one cycle. Therefore the total sum of cycles is $3$. 
    
    Let $f(x)=ax+b$ be an affine permutation of $\Z_{2^e}$. When $e=2$, we have $a\in\{1,3\}$. If $a=3$, then $\nu_2(a-1)\leq \nu_2(b)$ if and only if $b\in\{0,2\}$. In those cases, we have $\sum_{k=0}^2\frac{\varphi(p^k)}{\ord_{2^k}(3)}=3$ cycles in each map. If $b\in\{1,3\}$, each function generates $2$ cycles. When $a=1$, either $b=0$ and we have $2^2$ cycles, or $\nu_2(a-1)>\nu_2(b)$ and there are $2^{\nu_2(b)}$ cycles. Therefore, sum of the number of cycles is \begin{equation*}
        6+4+4+\sum_{b\in\Z_4\setminus\{0\}}2^{\nu_2(b)}=18.
    \end{equation*}

    For $e\geq 3$, we split the choices of $a$ into two cases: either $a=4j+3$ or $a=4j+1$, for some $j\in[0,2^{e-2}-1]$. In the first case, we have that $\nu_2(a-1)=1$. Therefore, $\nu_2(a-1)>\nu_2(b)$ if and only if $\nu_2(b)=0$ and there are $2^{e-1}$ elements $b\in\Z_{2^e}$ with that property. Moreover, from Lemma \ref{lif-exp} we obtain that $\nu_2(a^2-1)=3+\nu_2(j+1)\leq e+1$. We observe that for each $r\in[0,e-3]$ there are $2^{e-3-r}$ integers $j\in[0,2^{e-2}-1]$ with $\nu_2(j+1)=r$ and one such that $\nu_2(j+1)=e-2$. Hence, Theorem \ref{affinecycles} implies that the sum of the number of cycles generated by those maps is equal to \begin{equation}\label{2Part1}
        \sum_{j=0}^{2^{e-2}-1}\sum_{\substack{b\in\Z_{2^e}\\\nu_2(b)=0}}2^{\nu_p(j+1)+1}=2^{2e-2}+2^{e-1}\sum_{r=0}^{e-3}2^{e-2}=2^{2e-2}+(e-2)2^{2e-3}=e2^{2e-3}.
    \end{equation} 
    
    For any other $b\in\Z_{2^e}$, the sum of the number of cycles is given by \begin{align*}
        \sum_{j=0}^{2^{e-2}-1}\sum_{\substack{b\in\Z_{2^e}\\\nu_2(b)\geq1}}\sum_{k=0}^e\frac{\varphi(2^k)}{\ord_{2^k}(4j+3)}&=\sum_{j=0}^{2^{e-2}-1}2^{e-1}\left(2+\sum_{k=2}^e\frac{2^{k-1}}{2^{\max\{1,k-2-\nu_2(j+1)\}}}\right)\nonumber\\
        &=2^{2e-2}+2^{e-1}\left(\sum_{j=0}^{2^{e-2}-1}\sum_{k=0}^{e-2}2^{\min\{k,\nu_2(j+1)+1\}}\right).
    \end{align*} For each $r\in[0,e-3]$, there are $2^{e-3-r}$ integers $j\in[0,2^{e-2}-1]$ such that $\nu_2(j+1)=r$ and $\nu_2(j+1)=e-2$ if and only if $j=2^{e-2}-1$. Thus, we get \begin{align}
        2^{2e-2}+2^{e-1}\sum_{k=0}^{e-2}2^{\min\{k,e-1\}}+2^{e-1}\sum_{r=1}^{e-2}\sum_{k=0}^{e-2}2^{e-2-r+\min\{k,r\}}\nonumber\\
        &\hspace{-5.5cm}=2^{2e-2}+2^{e-1}(2^{e-1}-1)+2^{e-1}\sum_{r=1}^{e-2}\left(\sum_{k=0}^{r-1}2^{e-2-r+k}+\sum_{k=r}^{e-2}2^{e-2}\right)\nonumber\\
        &\hspace{-5.5cm}=2^{2e-1}-2^{e-1}+2^{e-1}\sum_{r=1}^{e-2}\left(2^{e-2}-2^{e-2-r}+(e-1-r)2^{e-2}\right)\nonumber\\
        &\hspace{-5.5cm}=2^{2e-1}-2^{e-1}+(e^2-e-2)2^{2e-4}-2^{e-1}(2^{e-2}-1)\nonumber\\
        &\hspace{-5.5cm}=(e^2-e+4)2^{2e-4}.\label{2Part3}
    \end{align}

    When $a=4j+1$ and $j\neq0$, we have that $\nu_2(a-1)=\nu_2(j)+2$. Hence, the sum of the number of cycles generated by those functions is given by \begin{align*}
        \sum_{j=1}^{2^{e-2}-1}\sum_{\substack{b\in\Z_{2^e}\\\nu_2(b)<\nu_2(j)+2}}2^{\nu_2(b)}+ \sum_{j=1}^{2^{e-2}-1}\sum_{\substack{b\in\Z_{2^e}\\\nu_2(b)\geq\nu_2(j)+2}}\sum_{k=0}^e\frac{\varphi(2^k)}{\ord_{2^k}(4j+1)}.
    \end{align*} 

    In the first summand, instead of summing over $b$, we sum over $t=\nu_2(b)$. Thus, using Lemma \ref{sumnup} we obtain
    \begin{equation}
        \sum_{j=1}^{2^{e-2}-1}\sum_{t=0}^{\nu_2(j)+1}2^{e-1}=2^{e-1}\sum_{j=1}^{2^{e-2}-1}(\nu_2(j)+2)=3\cdot2^{2e-3}-(e+1)2^{e-1}\label{2Part4a}
    \end{equation}
    In the second summand, we observe that the term is independent of $b$ and, for each $r$, there are $2^{e-r-2}$ choices of $b\in\Z_{2^e}$ such that $\nu_2(b)\geq r+2$. Summing over $r=\nu_2(j)$ and using Lemma \ref{ord} we get  \begin{align}
    \sum_{r=0}^{e-3}2^{2e-5-2r}\left(1+\sum_{k=1}^e\frac{2^{k-1}}{2^{\max\{0,k-2-r\}}}\right)\nonumber\\
        &\hspace{-2cm}=\frac{2^{2e-3}-2}{3}+\sum_{r=1}^{e-2}\left(\sum_{k=0}^{r}2^{2e-3-2r+k}+\sum_{k=r+1}^{e-1}2^{2e-3-r}\right)\nonumber\\
        &\hspace{-2cm}=\frac{2^{2e-3}-2}{3}+\sum_{r=1}^{e-2}2^{2e-2-r}-2^{2e-3-2r}+(e-r-1)2^{2e-3-r}\nonumber\\
        &\hspace{-2cm}=(e-1)2^{2e-3}-2^{e-1}\label{2Part4b}
    \end{align}
    When $j=0$, that is, $a=1$, we have that $\nu_2(b)\geq\nu_2(a-1)$ if and only if $b=0$. In that case we have the identity function and $2^{e}$ cycles. In any other case we have $2^{\nu_2(b)}$ cycles. Hence, we have \begin{equation}
        2^e+\sum_{b\in\Z_{2^e}\setminus\{0\}}2^{\nu_2(b)}=2^e+\sum_{t=0}^{e-1}2^{e-t-1}2^{e-1}=2^e+e2^{e-1}=(e+2)2^{e-1}\label{2Part5}
    \end{equation}

    Adding (\ref{2Part1}), (\ref{2Part3}), (\ref{2Part4a}), (\ref{2Part4b}) and (\ref{2Part5}) we obtain that the sum of the number of cycles of affine permutations over $\Z_{2^e}$ is $$(e^2+3e+8)2^{2e-4}.$$    
\end{proof}

\begin{theorem}\label{expeccycnumber2}
   The expectation for the number of cycles of a map $f\in\Ap(\Z_{2^e}$ selected uniformly at random is $$\Ex[\N(f,\Z_{2^e})\mid f\in\Ap(\Z_{2^e})]=\frac{(e^2+3e+8)}{8}.$$ 
\end{theorem}
\begin{proof}
    The total number of affine permutations in $\Z_{2^e}$ is $\varphi(2^e)\cdot 2^e=2^{2e-1}$ and, from Theorem \ref{cyc2e}, the total number of cycles generated by affine permutations is $(e^2+3e+8)2^{2e-4}$. Thus, the expectation for the number of cycles is given by $$\frac{1}{|\Ap(\Z_{2^e})|}\sum_{f\in\Ap(\Z_{2^e})}\N(f,\Z_{2^e})=\frac{(e^2+3e+8)2^{2e-4}}{2^{2e-1}}=\frac{e^2+3e+8}{8}.$$
\end{proof}

\begin{theorem}\label{secmom2e}
    $$\sum_{f\in\Ap(\Z_{2^e})}\sum_{\ell\in\L(f,\Z_{2^e})}\ell^2\N_\ell(f,\Z_{2^e})=\frac{2^{2e}}{18}(2^{2e+2}+6e+5).$$
\end{theorem}
\begin{proof}
    For any prime power $p^e$ and affine permutation $f(x)=ax+b$ over $\Z_{p^e}$, if $\nu_p(a-1)\leq \nu_p(b)$, Theorem \ref{altaffinecycles} implies that $$\sum_{\ell\in\L(f,\Z_{2^e})}\ell^2\N_\ell(f,\Z_{2^e})=\sum_{k=0}^e\varphi(p^k)\ord_{p^k}(a).$$ If $\nu_p(a-1)>\nu_p(b)$, Theorem \ref{affinecycles} implies that $$\sum_{\ell\in\L(f,\Z_{2^e})}\ell^2\N_\ell(f,\Z_{2^e})=p^{e+s},$$ where $s$ depends on $\nu_p(b)$ and $\nu_p(a^2-1)$. We now follow the same structure as the proof of Theorem \ref{cyc2e}. 

    When $p=2$, if $e=1$ the only affine permutations are the identity and $x\mapsto x+1$. Thus $$\sum_{f\in\Ap(\Z_{2})}\sum_{\ell\in\L(f,\Z_{2})}\ell^2\N_\ell(f,\Z_{2})=1\cdot2+4\cdot1=6.$$
    
    When $e=2$, we have $a\in\{1,3\}$. If $a=3$, then $\nu_2(a-1)\leq \nu_2(b)$ if and only if $b\in\{0,2\}$. In each of those cases, we have $\sum_{k=0}^2{\varphi(p^k)}{\ord_{2^k}(3)}=6$. If $b\in\{1,3\}$, we have $p^{e+s}=2^3=8$. When $a=1$, either $b=0$ and we have $2^2$ cycles of length $1$, or $\nu_2(a-1)>\nu_2(b)$ and $2^{e+s}=2^{2e-\nu_2(b)}$. Therefore, we obtain \begin{equation*}
        \sum_{f\in\Ap(\Z_{2^2})}\sum_{\ell\in\L(f,\Z_{2^2})}\ell^2\N_\ell(f,\Z_{2^2})=12+16+4+\sum_{b\in\Z_4\setminus\{0\}}2^{4-\nu_2(b)}=72.
    \end{equation*}
    
    For $e\geq3$, if $a=4j+3$ and $\nu_2(b)=0$, for some $j\in[0,2^{e-2}-1]$, then $r=\nu_2(j+1)$ and we have $2^{r+1}$ cycles of length $2^{e-r-1}$. Thus, $$\sum_{\ell\in\L(f,\Z_{2^e})}\ell^2\N_\ell(f,\Z_{2^e})=2^{2e-r-1}.$$ 
    
    There are $2^{e-1}$ choices of $b$ such that $\nu_2(b)=0$ and $2^{e-3-r}$ choices of $j$ with $\nu_2(j+1)=r$, for each $r\in[0,e-3]$, and one with $\nu_2(j+1)=e-2$. Hence, we have 
    \begin{align}
   \sum_{j=0}^{2^{e-2}-1}\sum_{\substack{b\in\Z_{2^e}\\\nu_2(b)=0}}2^{2e-\nu_2(j+1)-1}=2^{e-1}\left(2^{2e-(e-2)-1}+\sum_{r=0}^{e-3}2^{e-3-r}2^{2e-r-1}\right)=\frac{2^{4e-3}+2^{2e}}{3}.\label{S1}
    \end{align}

    If $a=4j+3$ and $\nu_2(b)>0$, then \begin{align*}
        \sum_{j=0}^{2^{e-2}-1}\sum_{\substack{b\in\Z_{2^e}\\\nu_2(b)\geq1}}\sum_{k=0}^e\varphi(p^k)\ord\nolimits_{p^k}(4j+3)=&
        2^{e-1}\sum_{r=0}^{e-3}2^{e-3-r}\left(2+\sum_{k=2}^e2^{k-1}2^{\max\{1,k-2-r\}}\right)\nonumber\\
        &+2^{e-1}\left(2+\sum_{k=2}^e2^{k}\right).
    \end{align*}
    The second term comes from the unique $j$ such that $\nu_2(j+1)=e-2$, and it is equal to $2^{2e}-2^{e}$. In the first term, we have that \begin{align*}
        2+\sum_{k=2}^e2^{k-1}2^{\max\{1,k-2-r\}}&=2+\sum_{k=2}^{r+3}2^{k}+\sum_{k=r+4}^e2^{2k-r-3}\\&=2+2^{r+4}-4+\frac{2^{2e-r-1}-2^{r+5}}{3}\\&=\frac{2^{2e-r-1}+2^{r+4}-6}{3}.
    \end{align*}
    Multiplying by $2^{2e-4-r}$ and summing over $r$ we get \begin{align*}
        \sum_{r=0}^{e-3}2^{2e-4-r}\frac{2^{2e-r-1}+2^{r+4}-6}{3}&=\sum_{r=0}^{e-3}\frac{2^{4e-2r-5}+2^{2e}-3\cdot 2^{2e-3-r}}{3}\\
        &=\frac{1}{3}\left(\frac{2^{4e-3}-2^{2e+1}}{3}+(e-2)2^{2e}-3(2^{2e-2}-2^{e})\right).
    \end{align*} Adding both terms, we get \begin{align*}
    \frac{1}{9}\left(9(2^{2e}-2^{e})+2^{4e-3}-2^{2e+1}+3(e-2)2^{2e}-9(2^{2e-2}-2^{e})\right), 
    \end{align*} which simplifies to \begin{equation}
        \frac{2^{2e-3}}{9}\left(2^{2e}+24e-10\right)\label{S2}.
    \end{equation}

    For $a=4j+1$ and $j\neq 0$, let $r=\nu_2(j)$. If $\nu_2(b)=t<r+2$, then $2^{e+s}=2^{2e-t}$. There are $2^{e-t-1}$ elements $b$ with $\nu_2(b)=t$, for $t\in[0,e-1]$. Moreover, there are $2^{e-r-3}$ integers $j\in[0,2^{e-2}-1]$ with $\nu_2(j)=r$, for $r\in[0,e-3]$. Thus, we have \begin{align}
    \sum_{r=0}^{e-3}2^{e-r-3}\sum_{t=0}^{r+1}2^{e-t-1}2^{2e-t}&=\sum_{r=0}^{e-3}2^{2e-r-4}\sum_{t=0}^{r+1}2^{2e-2t}\nonumber\\
    &=\sum_{r=0}^{e-3}2^{2e-r-4}\left(\frac{2^{2e+2}-2^{2e-2r-2}}{3}\right)\nonumber\\
    &=\frac{1}{3}\left(2^{4e-1}-2^{3e+1}-\frac{2^{4e-3}-2^{e+3}}{7}\right)\nonumber\\
    &=\frac{9}{7}2^{4e-3}-\frac{1}{3} 2^{3e+1}+\frac1{21}2^{e+3}.\label{S3}
    \end{align}

    If $\nu_2(b)=t\geq r+2$, we have \begin{align*}
        \sum_{k=0}^e\varphi(p^k)\ord\nolimits_{p^k}(4j+1)&=2+\sum_{k=2}^e2^{k-1+\max\{0,k-2-r\}}\\
        &=2+\sum_{k=2}^{r+2}2^{k-1}+\sum_{k=r+3}^e2^{2k-r-3}\\
        &=\frac{2^{2e-r-1}+2^{r+2}}{3}.
    \end{align*}
    There are $2^{e-r-2}$ choices for $b$ with $\nu_2(b)\geq r+2$ and, for every $r\in[0,e-3]$, there are $2^{e-r-3}$ integers $j\in[0,2^{e-2}-1]$ with $\nu_2(j)=r$. Hence, summing over $r\in[0,e-3]$ we have  \begin{align}
        \sum_{r=0}^{e-3}2^{e-r-3}2^{e-r-2}\left(\frac{2^{2e-r-1}+2^{r+2}}{3}\right)&=\frac{1}{3}\sum_{r=0}^{e-3}(2^{4e-3r-6}+2^{2e-r-3})\nonumber\\
        &=\frac{1}{3}\left(\frac{2^{4e-3}-2^{e+3}}{7}+2^{2e-2}-2^e\right)\nonumber\\
        &=\frac{1}{21}\left(2^{4e-3}+7\cdot 2^{2e-2}-15\cdot2^{e}\right).\label{S4}
    \end{align}

    When $a=1$ we have $\sum_{\ell\in\L(f,\Z_{2^e})}\ell^2\N_\ell(f,\Z_{2^e})=2^{e+s}=~2^{2e-\nu_2(b)}$ for every $b\in\Z_{2^e}\setminus\{0\}$ and $\sum_{\ell\in\L(f,\Z_{2^e})}\ell^2\N_\ell(f,\Z_{2^e})=2^e$ if $b=0$. Thus, we obtain \begin{align}
        2^e+\sum_{b\in\Z_{2^e}\setminus\{0\}}2^{2e-\nu_2(b)}&=2^e+\sum_{t=0}^{e-1}2^{e-t-1}2^{2e-t}\nonumber\\&=2^e+\sum_{t=0}^{e-1}2^{3e-2t-1}\nonumber\\&=\frac{2^{3e+1}+2^e}{3}.\label{S5}
    \end{align}

    Adding (\ref{S1}), (\ref{S2}), (\ref{S3}), (\ref{S4}) and (\ref{S5}), we conclude that $$\sum_{f\in\Ap(\Z_{2^e})}\sum_{\ell\in\L(f,\Z_{2^e})}\ell^2\N_\ell(f,\Z_{2^e})=\frac{2^{2e}}{18}(2^{2e+2}+6e+5).$$
\end{proof}

\begin{theorem}\label{expeccyclength2}
    The expectation for the cycle length of an element $x\in\Z_{2^e}$ over a function $f\in\Ap(\Z_{2^e})$, both selected uniformly at random, is $$\Ex[\ell_x(f,\Z_{2^e})]=\frac{2^{2e+2}+6e+5}{9\cdot 2^e}.$$ 
\end{theorem}
\begin{proof}
    The total number of affine permutations in $\Z_{2^e}$ is $\varphi(2^e)\cdot 2^e=2^{2e-1}$ and the total number of points is $2^e$. From the previous theorem we conclude that the expectation is $$\frac{2^{2e}(2^{2e+2}+6e+5)}{18\cdot 2^{3e-1}}=\frac{2^{2e+2}+6e+5}{9\cdot 2^e}.$$
\end{proof}

\subsection{Stochastic properties for odd primes}

\begin{theorem}\label{cycpe}
    Let $p$ be an odd prime and $e$ a positive integer. Then  $$\sum_{f\in\Ap(\Z_{p^e})}\N(f,\Z_{p^e})=p^{2e-2}(p-1)\left(p+e+\frac{e((p-1)e+p+1)}{2}F\right),$$ where $F=F(p-1)=\sum\limits_{\substack{d\mid p-1\\d\neq1}}\frac{\varphi(d)}{d}$. 
\end{theorem}
\begin{proof}
    Let $f\in\Ap(\Z_{p^e})$ be given by $f(x)=ax+b$. From Theorem \ref{affinecycles} we have $$\sum_{\ell\in \L(f,\Z_{p^e})}\N_\ell(f,\Z_{p^e})=\begin{cases}
        \sum_{k=0}^e\frac{\varphi(p^k)}{\ord_{p^k}(a)}&\text{if $\nu_p(a-1)\leq\nu_p(b)$,}\\
        p^{\nu_p(b)}&\text{otherwise}.
    \end{cases}$$
    We begin with the case where $e=1$. If $a\in[2,p-1]$, then $\nu_p(a-1)\leq \nu_p(b)$ for every $b\in\Z_p$. Thus 
    \begin{equation}\label{e1inot1}
        \sum_{b\in\Z_{p}}\sum_{i=2}^{p-1}\sum_{\ell\in \L(f,\Z_{p^e})}\N_\ell(f,\Z_{p^e})=\sum_{b\in\Z_{p}}\sum_{i=2}^{p-1}1+\frac{\varphi(p)}{\ord_{p}(i)}=p(p-2)+p(p-1)\sum_{i=2}^{p-1}\frac{1}{\ord_p(i)}.
    \end{equation}
    For each $d$ divisor of $p-1$, there are $\varphi(d)$ elements with order $d$ in $U(\Z_p)$. Let $F=\sum\limits_{\substack{d\mid p-1\\d\neq1}}\frac{\varphi(d)}{d}$. Thus, (\ref{e1inot1}) is equal to \begin{equation}
        p(p-2)+p(p-1)F\label{e1Part1}
    \end{equation}
    If $a=1$, then $\nu_p(a-1)\leq \nu_p(b)$ if and only if $b=0$. Thus \begin{equation}
        \sum_{b\in\Z_{p}}\sum_{\ell\in \L(f,\Z_{p^e})}\N_\ell(f,\Z_{p^e})=1+(p-1)+\sum_{b\in\Z_p\setminus\{0\}}p^{\nu_p(b)}=2p-1.\label{e1Part2}
    \end{equation}
    Adding (\ref{e1Part1}) and (\ref{e1Part2}) we get $$(p-1)(pF+p+1).$$

    For $e>1$, we notice that any $a\in\U(\Z_{p^e})$ can be written as $a=jp+i$ with $j\in[0,p^{e-1}-1]$ and $i\in[1,p-1]$. When $i\neq 1$, we have that $\nu_p(a-1)\leq\nu_p(b)$ for every $b\in\Z_{p^e}$. Therefore, it follows that 
    \begin{equation}\label{inot1}
        \sum_{b\in\Z_{p^e}}\sum_{j=0}^{p^{e-1}-1}\sum_{i=2}^{p-1}\sum_{\ell\in \L(f,\Z_{p^e})}\N_\ell(f,\Z_{p^e})=\sum_{b\in\Z_{p^e}}\sum_{j=0}^{p^{e-1}-1}\sum_{i=2}^{p-1}\sum_{k=0}^e\frac{\varphi(p^k)}{\ord_{p^k}(jp+i)}.
    \end{equation}
   Let $d_i=\ord_p(i)$, for $i\in[2,p-1]$. We observe that $\ord_p(jp+i)=\ord_p(i)$. For a fixed $i$, let $t_j=\nu_p((jp+i)^{d_i}-1)$. When $k=0$, we have that $\frac{\phi(p^k)}{\ord_{p^k}(jp+i)}=1$ and $\sum_{j=0}^{p^{e-1}-1}1=p^{e-1}$. Otherwise, for any $k\in[1,e]$, we have that $\ord_{p^k}(jp+i)=d_ip^{\max\{0,k-t_j\}}$ and $$\frac{\varphi(p^k)}{\ord_{p^k}(jp+i)}=\frac{p^{k}-p^{k-1}}{d_ip^{\max\{0,k-t_j\}}}=\frac{(p-1)p^{\min\{k-1,t_j-1\}}}{d_i}.$$
    As $t_j\geq 1$ for every $j$ and $d_i<p$, Lemma \ref{hensel} implies that, for every $r\in[1,e]$, there are $p^{e-r}$ choices of $j\in[0,p^{e-1}-1]$ such that $t_j\geq r$. Thus, there are exactly $p^{e-r}-p^{e-r-1}$ choices such that $t_j=r$, for $r\in[1,e-1]$ and $1$ such that $t_j\geq e$. Hence \begin{align*}
        \sum_{j=0}^{p^{e-1}-1}p^{\min\{k-1,t_j-1\}}&=p^{e-k}p^{k-1}+\sum_{r=1}^{k-1}p^{r-1}(p^{e-r}-p^{e-r-1})\\
        &=p^{e-1}+(k-1)(p^{e-1}-p^{e-2})\\
        &=kp^{e-1}-(k-1)p^{e-2}.
    \end{align*}   
    Therefore, we have $$\sum_{k=1}^{e}kp^{e-1}-(k-1)p^{e-2}=\frac{e(e+1)}{2}p^{e-1}-\frac{e(e-1)}{2}p^{e-2}$$ and (\ref{inot1}) is equal to $$\sum_{b\in\Z_{p^e}}\sum_{i=2}^{p-1}\left(p^{e-1}+\frac{p-1}{d_i}\left(\frac{e(e+1)}{2}p^{e-1}-\frac{e(e-1)}{2}p^{e-2}\right)\right).$$

    To sum over $i$, we observe that for every $d$ divisor of $p-1$, there are $\varphi(d)$ choices of $i\in[1,p-1]$ such that $d_i=d$. Thus,    $$\sum_{i=2}^{p-1}\frac{1}{d_i}=\sum_{\substack{d\mid p-1\\d>1}}\frac{\varphi(d)}{d}.$$ Also summing over $b\in\Z_{p^e}$, this concludes that (\ref{inot1}) is equal to \begin{equation}
        p^{2e-1}(p-2)+p^e(p-1)\left(\frac{e(e+1)}{2}p^{e-1}-\frac{e(e-1)}{2}p^{e-2}\right)\sum_{\substack{d\mid p-1\\d>1}}\frac{\varphi(d)}{d}.\label{part1}
    \end{equation}

    Now, we look at the case where $i=1$ and $j>0$. We have that $\nu_p(a-1)\leq\nu_p(b)$ if and only if $\nu_p(b)> \nu_p(j)$. Thus, we have     \begin{align}\label{i1}
        \sum_{b\in\Z_{p^e}}\sum_{j=1}^{p^{e-1}-1}\sum_{\ell\in \L(f,\Z_{p^e})}\N_\ell(f,\Z_{p^e})=&\sum_{j=1}^{p^{e-1}-1}\sum_{\substack{b\in\Z_{p^e}\\\nu_p(b)> \nu_p(j)}}\sum_{k=0}^e\frac{\varphi(p^k)}{\ord_{p^k}(jp+1)}\nonumber \\  &+ \sum_{j=1}^{p^{e-1}-1}\sum_{\substack{b\in\Z_{p^e}\\\nu_p(b)\leq \nu_p(j)}}p^{\nu_p(b)}.
    \end{align}

    In the first summand, when $k=0$ we have $\frac{\varphi(p^k)}{\ord_{p^k}(jp+1)}=1$. Otherwise, we get $\frac{\varphi(p^k)}{\ord_{p^k}(jp+1)}=(p-1)p^{\min\{k-1,\nu_p(j)\}}$. For each $r\in[0,e-2]$, there are $p^{e-r-2}(p-1)$ choices of $j\in[1,p^{e-1}-1]$ with $\nu_p(j)=r$. Moreover, there are $p^{e-r-1}$ choices for $b\in\Z_{p^e}$ such that $\nu_p(b)> r$. Therefore, we have that \begin{align}
        \sum_{j=1}^{p^{e-1}-1}\sum_{\substack{b\in\Z_{p^e}\\\nu_p(b)> \nu_p(j)}}\sum_{k=0}^e\frac{\varphi(p^k)}{\ord_{p^k}(jp+1)}\nonumber\\
        &\hspace{-2cm}=\sum_{r=0}^{e-2}p^{2e-2r-3}(p-1)\left(1+\sum_{k=1}^{r+1}(p-1)p^{k-1}+\sum_{k=r+2}^{e}(p-1)p^{r}\right)\nonumber\\
        &\hspace{-2cm}=\sum_{r=0}^{e-2}p^{2e-2r-3}(p-1)(1+(p^{r+1}-1)+(e-r-1)(p-1)p^r)\nonumber\\
        &\hspace{-2cm}=p^{e-2}\sum_{r=0}^{e-2}p^{e-r-1}(p-1)(p+(e-r-1)(p-1))\nonumber\\
        &\hspace{-2cm}=p^{e-1}(p-1)(ep^{e-1}-1)\label{Part2a}
    \end{align}

    In the second summand of (\ref{i1}), using Lemma \ref{sumnup} we have that \begin{align}
        \sum_{j=1}^{p^{e-1}-1}\sum_{\substack{b\in\Z_{p^e}\nonumber\\\nu_p(b)\leq \nu_p(j)}}p^{\nu_p(b)}&=\sum_{j=1}^{p^{e-1}-1}\sum_{r=0}^{\nu_p(j)}p^{r}(p^{e-r}-p^{e-1-r})\nonumber\\&=(p^e-p^{e-1})\sum_{j=1}^{p^{e-1}-1}\nu_p(j)+1\nonumber\\&=(p^e-p^{e-1})(p^{e-1}-1+\frac{p^{e-1}-1}{p-1}-(e-1))\nonumber\\&=p^e(p^{e-1}-1)-p^{e-1}(p-1)(e-1).\label{Part2b}
    \end{align}

    Now when $i=1$ and $j=0$, we have that $\nu_p(b)\geq \nu_p(a-1)$ if and only if $b=0$. Hence, we have \begin{align}
        \sum_{b\in\Z_{p^e}}\sum_{\ell\in \L(f,\Z_{p^e})}\N_\ell(f,\Z_{p^e})&=\sum_{k=0}^e\frac{\varphi(p^k)}{\ord_{p^k}(1)} + \sum_{b\in\Z_{p^e}\setminus\{0\}}p^{\nu_p(b)}\nonumber\\&=p^e+\sum_{r=0}^{e-1}p^{r}(p^{e-r}-p^{e-r-1})\nonumber\\&=p^e+e(p^e-p^{e-1}).\label{Part3}
    \end{align}

    Adding (\ref{part1}), (\ref{Part2a}), (\ref{Part2b}) and (\ref{Part3}), we obtain \begin{align*}
      p^{2e-2}(p-1)(p+e)+p^e(p-1)\left(\frac{e(e+1)}{2}p^{e-1}-\frac{e(e-1)}{2}p^{e-2}\right)\sum_{\substack{d\mid p-1\\d>1}}\frac{\varphi(d)}{d}\\=p^{2e-2}(p-1)\left(p+e+\frac{e((p-1)e+p+1)}{2}F\right)
    \end{align*}
\end{proof}

\begin{theorem}\label{expeccycnumberp}
   Let $p$ be an odd prime and $e$ a positive integer. The expectation for the number of cycles in the functional graph of a map $f\in\Ap(\Z_{p^e})$ selected uniformly at random is $$\Ex[\N(f,\Z_{p^e})]=1+\frac{e}{p}+\frac{e((p-1)e+p+1)}{2p}F.$$
\end{theorem}
\begin{proof}
    This result follows by dividing the total number of cycles from Theorem \ref{cycpe} by $\varphi(p^e)\cdot p^e=p^{2e-1}(p-1)$, the total number of affine permutations in $\Z_{p^e}$.
\end{proof}

\begin{theorem}\label{secmompe}
     Let $p$ be an odd prime and $e$ a positive integer. Then  \begin{align*}
         \sum_{f\in\Ap(\Z_{p^e})}\sum_{\ell\in\L(f,\Z_{p^e})}\ell^2\N_\ell(f,\Z_{p^e})\\&\hspace{-5cm}
         =\frac{p^{2e-1}}{(p+1)^2}\left((p+1)p^{2e+1}+p^3-2p-1+G(p^{2e+1}+e(p^2-1)-p)\right),
     \end{align*} where $G=G(p-1)=\sum\limits_{\substack{d\mid p-1\\d\neq 1}}d\cdot\varphi(d)$. 
\end{theorem}
\begin{proof}
    For an affine permutation $f(x)=ax+b$ over $\Z_{p^e}$, we have that  $$\sum_{\ell\in\L(f,\Z_{p^e})}\ell^2\N_\ell(f,\Z_{p^e})=\sum_{k=0}^e\varphi(p^k)\ord\nolimits_{p^k}(a), $$ if $\nu_p(a-1)\leq \nu_p(b)$, and $$\sum_{\ell\in\L(f,\Z_{p^e})}\ell^2\N_\ell(f,\Z_{p^e})=p^{2e-\nu_p(b)},$$ if $\nu_p(a-1)>\nu_p(b)$.

    Suppose $e=1$. If $a\in[2,p-1]$, then $\nu_p(a-1)\leq\nu_p(b)$ for all $b\in\Z_p$ and $$\sum_{\ell\in\L(f,\Z_{p})}\ell^2\N_\ell(f,\Z_{p})=1+(p-1)\ord\nolimits_p (a).$$ Summing over $b$ and $a$ we have $$\sum_{a\in[2,p-1]}\sum_{b\in\Z_p}1+(p-1)\ord\nolimits_p(a)=p(p-2)+p(p-1)\sum_{\substack{d\mid p-1\\d\neq 1}}d\varphi(d).$$ Defining $G=\sum\limits_{\substack{d\mid p-1\\d\neq1}}d\varphi(d)$, we have \begin{equation}
        p(p-2)+p(p-1)G.\label{e1a}
    \end{equation}

    If $a=1$, when $b=0$ we have the identity map and $\sum_{\ell\in\L(f,\Z_{p})}\ell^2\N_\ell(f,\Z_{p})=p.$ Otherwise, we have $\sum_{\ell\in\L(f,\Z_{p})}\ell^2\N_\ell(f,\Z_{p})=p^2$ for every $b\in\Z_p\setminus\{0\}$, thus \begin{equation}
        \sum_{b\in\Z_p}\sum_{\ell\in\L(f,\Z_{p^e})}\ell^2\N_\ell(f,\Z_{p^e})=p+(p-1)p^2.\label{e1b}
    \end{equation} Adding (\ref{e1a}) and (\ref{e1b}) we get $$p(p-1)(G+p+1).$$

    Now, suppose $e>1$. If $a=jp+i$, where  $j\in[0,p^{e-1}-1]$ and $i\in[2,p-1]$, then $\nu_p(a-1)\leq\nu_p(b)$ for every $b\in\Z_{p^e}$. Let $d_i=\ord_p(i)$ and $t_j=\nu_p((jp+i)^{d_i}-1)$. When $k=0$, we have \begin{equation}
        \sum_{b\in\Z_{p^e}}\sum_{j=0}^{p^{e-1}-1}\sum_{i=2}^{p-1}\varphi(p^k)\ord\nolimits_{p^k}(jp+i)=p^{2e-1}(p-2).\label{E1}
    \end{equation}
    For $k=1$ we have \begin{equation}
        \sum_{b\in\Z_{p^e}}\sum_{j=0}^{p^{e-1}-1}\sum_{i=2}^{p-1}\varphi(p^k)\ord\nolimits_{p^k}(jp+i)=p^{2e-1}(p-1)G.\label{E2}
    \end{equation}
    And, for $k\in[2,e]$, we have $\varphi(p^k)\ord\nolimits_{p^k}(jp+i)=(p-1)p^{\max\{k-1,2k-t_j-1\}}d_i$. For every $r\in[1,e-1]$, there are $p^{e-r}-p^{e-r-1}$ values for $j$ such that $t_j=r$ and one value with $t_j\geq e$. For a fixed $k$, we have \begin{align*}
        \sum_{j=0}^{p^{e-1}-1}d_i(p-1)p^{\max\{k-1,2k-t_j-1\}}&=d_i(p-1)\left(p^{e-1}+\sum_{r=1}^{k-1}(p^{e-r}-p^{e-r-1})p^{2k-r-1}\right)\\
        &=d_i(p-1)\left(p^{e-1}+(p-1)\sum_{r=1}^{k-1}p^{e+2k-2r-2}\right)\\
        &=d_i(p-1)\left(p^{e-1}+\frac{p^{e+2k-2}-p^{e}}{p+1}\right)\\
        &=d_i(p-1)\frac{p^{e+2k-2}+p^{e-1}}{p+1}.
    \end{align*}
    Summing over $k\in[2,e]$, we get \begin{align*}
        \sum_{k=2}^ed_i(p-1)\frac{p^{e+2k-2}+p^{e-1}}{p+1}&=\frac{d_ip^{e-1}}{(p+1)^2}\left(p^{2e+1}-p^3+(e-1)(p^2-1)\right).
    \end{align*} And summing over $i\in[2,p-1]$ and $b\in\Z_{p^e}$, we obtain \begin{align}
        \sum_{b\in\Z_{p^e}}\sum_{i=2}^{p-1}\frac{d_ip^{e-1}}{(p+1)^2}\left(p^{2e+1}-p^3+(e-1)(p^2-1)\right)\nonumber\\&\hspace{-5cm}=\frac{p^{2e-1}}{(p+1)^2}\left(p^{2e+1}-p^3+(e-1)(p^2-1)\right)G.\label{E3}
    \end{align} Adding (\ref{E1}), (\ref{E2}) and (\ref{E3}) we have \begin{equation}
        p^{2e-1}\left(p-2+\frac{p^{2e+1}+e(p^2-1)-p}{(p+1)^2}G\right).\label{Eq1}
    \end{equation}

    If $i=1$ and $j>0$, let $r=\nu_p(j)$. We have that $\nu_p(a-1)=r+1$ and either $r+1\leq \nu_p(b)$ or $r\geq\nu_p(b)$. In the first case we have \begin{align}
        \sum_{\substack{b\in\Z_{p^e}\\\nu_p(b)>\nu_p(j)}}\sum_{j=1}^{p^{e-1}-1}\sum_{k=0}^e\varphi(p^k)\ord\nolimits_{p^k}(jp+1)\nonumber\\
        &\hspace{-4cm}=\sum_{r=0}^{e-2}p^{2e-2r-3}(p-1)\left(1+\sum_{k=1}^{r+1}(p-1)p^{k-1}+\sum_{k=r+2}^e(p-1)p^{2k-r-2}\right)\nonumber\\
        &\hspace{-4cm}=\sum_{r=0}^{e-2}p^{2e-2r-3}(p-1)\left(\frac{p^{r+1}+p^{2e-r}}{p+1}\right)\nonumber\\
        &\hspace{-4cm}=(p-1)\sum_{r=0}^{e-2}\left(\frac{p^{2e-r-2}+p^{4e-3r-3}}{p+1}\right)\nonumber\\
        &\hspace{-4cm}=\frac{1}{p+1}\left(p^{2e-1}-p^e+\frac{p^{4e}-p^{e+3}}{p^2+p+1}\right)\label{Eq2}
    \end{align} In the second case, let $t=\nu_p(b)$. We have that \begin{align}
        \sum_{j=1}^{p^{e-1}-1}\sum_{\substack{b\in\Z_{p^e}\\\nu_p(b)\leq\nu_p(j)}}p^{2e-\nu_p(b)}
        &=\sum_{r=0}^{e-2}(p-1)p^{e-r-2}\sum_{t=0}^r(p-1)p^{e-t-1}p^{2e-t}\nonumber\\
        &=(p-1)^2\sum_{r=0}^{e-2}\sum_{t=0}^rp^{4e-r-2t-3}\nonumber\\
        &=\frac{p^{4e}-p^{3e+1}}{p+1}-\frac{p^{4e-p^{e+3}}}{(p+1)(p^2+p+1)}\label{Eq3}.
    \end{align} 

    If $i=1$ and $j=0$, when $b=0$ we have the identity map and \begin{equation*}
        \sum_{\ell\in\L(f,\Z_{p^e})}\ell^2\N_\ell(f,\Z_{p^e})=p^e.
    \end{equation*} Adding the other values of $b\in\Z_{p^e}\setminus\{0\}$ we have \begin{align}
        p^e+\sum_{b\in\Z_{p^e}\setminus\{0\}}p^{2e-\nu_p(b)}
        =p^e+\sum_{t=0}^{e-1}(p-1)p^{e-t-1}p^{2e-t}=\frac{p^{3e+1}+p^e}{p+1}.\label{Eq4}
    \end{align} 
    Adding (\ref{Eq1}), (\ref{Eq2}), (\ref{Eq3}) and (\ref{Eq4}), the result follows.    
\end{proof}

\begin{theorem}\label{expeccyclengthp}
   Let $p$ be an odd prime and $e$ a positive integer. The expectation for cycle length of an element $x\in\Z_{p^e}$ in the functional graph of a map $f\in\Ap(\Z_{p^e})$, both selected uniformly at random, is $$\Ex[\ell_x(f,\Z_{p^e})]=
   \frac{p(p+1+G)}{(p-1)(p+1)^2}p^e+\frac{eG}{(p+1)p^e}+\frac{p^3-2p-1-pG}{(p-1)(p+1)^2p^e}$$
\end{theorem}
\begin{proof}
    This result follows by dividing the total number of cycles from Theorem \ref{secmompe} by $p^e\cdot |\Ap(\Z_{p^e})|=p^{3e-1}(p-1)$.
\end{proof}

\begin{remark}
    Lemma \ref{Theorem2} determines the functional graph of an affine map $g$ over $\Z_n$, for $n$ composite, as the product of the functional graph of $g$ over $\Z_{p_i^{e_i}}$, for every $p_i^{e_i}$ in the prime decomposition of $n$. However, from Lemma \ref{lcm}, we have that the total number of cycles depends on the $\gcd$s of each component's cycle length. Therefore, a closed formula for the total number of cycles generated by all affine maps over $\Z_n$ would rely on the sum of the $\gcd$s over all combinations of cycle lengths of each component, making it substantially more complex.
\end{remark}

\section{A comparison with random permutations}\label{Random}
Let $\sigma$ be a permutation of the set of $n$ integers $[1,n]$ and $x\in[1,n]$, both selected uniformly at random. Let $\N(\sigma,n)$ be the number of cycles of $\sigma$ and $\ell_x(\sigma,n)$ be the length of the cycle that contains $x$ in $\sigma$. A classical result from random permutation theory is the expectation for those values (see \cite[Chapter 1.3]{Stanley_2011}), which is given by \begin{equation*}
    \Ex[\N(\sigma,n)]=H_n=\sum_{k=1}^n\frac{1}{k}\quad\text{and}\quad\Ex[\ell_x(\sigma,n)]=\frac{n+1}{2}.
\end{equation*}

We first observe the ratio between the expectation for affine permutations over $\Z_{p^e}$ and random permutations. Let $$R^{\Ap}_\N(p^e)=\frac{ \Ex[\N(g,\Z_{p^e})]}{ \Ex[\N(\sigma,p^e)]}\quad \text{and}\quad  R^{\Ap}_\ell(p^e)=\frac{\Ex[\ell_t(g,\Z_{p^e})]}{\Ex[\ell_x(\sigma,p^e)]}.$$  Using the results of Section \ref{SAff}, we first investigate their asymptotic behavior when $e$ increases and $p$ is fixed.
\begin{theorem}
Let $p$ be a fixed prime. As $e$ tends to infinity, the ratios $R^{\Ap}_\N(p^e)$ and $R^{\Ap}_\ell(p^e)$ satisfy    
$$R_\N^{\Ap}(2^e)=\frac{e}{8\log2}+\O(1),\quad\text{and}\quad R_\ell^{\Ap}(2^e)=\frac{8}{9}+\O(2^{-e}),
$$
when $p=2$, and
$$R_\N^{\Ap}(p^e)=\frac{(p-1)F}{2p\log p}e+\O(1),\quad\text{and}\quad R_\ell^{\Ap}(p^e)=\frac{2p(p+1+G)}{(p-1)(p+1)^2}+\O(ep^{-e}),$$
when $p$ is odd, where $F=\sum_{\substack{r\mid p-1\\r>1}}\frac{\varphi(r)}r$ and $G=\sum_{\substack{r\mid p-1\\r>1}}r\varphi(r)$ and the implied constants depend only on $p$.
\end{theorem}
\begin{proof}
    The result follows directly from the asymptotic expansion of the harmonic numbers and Theorems \ref{expeccycnumber2}, \ref{expeccyclength2}, \ref{expeccycnumberp} and \ref{expeccyclengthp}.
\end{proof}

In particular, the ratio $R_\N^{\Ap}(p^e)$ grows linearly with $e$, while $R_\ell^{\Ap}(p^e)$ tends to a constant determined by $p$.

The analysis when $p$ grows depends heavily on the prime factorization of $p-1$ and cannot be obtained as simply in general. We provide an asymptotic result for \textit{safe primes}. A prime number $p$ is a safe prime if it can be written as $p=2m+1$, for a prime $m$. Although their infinitude remains an open problem, numerical evidence and the Hardy–Littlewood conjectures point towards the existence of infinitely many safe primes.

\begin{theorem}
    Let $e$ be a fixed positive integer. For every safe prime $p=2r+1$, where $r$ is an odd prime, we have
    $$R_\N^{\Ap}(p^e)=\frac{1+e(e+1)}{e\log p}+\O((\log p)^{-2})\quad\text{and}\quad R_\ell^{\Ap}(p^e)=\frac{3}{2}+\O(p^{-1}),$$
    as $p$ increases through safe primes, where the implied constants depend on $e$.
\end{theorem}

\begin{proof}
    We observe that if $p=2m+1$, for some odd prime $m$, then $$F=\sum_{r\mid p-1\\r\neq1}\frac{\varphi(r)}{r}=2-\frac{3}{2m}=2-\frac{3}{p-1}$$ and $$G=\sum_{r\mid p-1\\r\neq1}r\varphi(r)=3m^2-3m+2=\frac{3p^2-12p+17}4.$$ The result follows from Theorems \ref{expeccycnumberp}, \ref{expeccyclengthp}.
\end{proof}

Consequently, along any unbounded sequence of safe primes, the ratio between the expected number of cycles of affine permutations and random permutations becomes negligible, while the ratio of the expected length of the cycle containing a random element approaches $\frac{3}{2}$.

For generalized cyclotomic permutations, we observe that the coset permutation built in their dynamics contributes exactly with factors equal to the expectations for random permutations in the expected number of cycles and the expected cycle length of a random point. Therefore, the ratio between those two properties for GCPs and random permutations is given by \begin{equation*}
    R^{GCP}_\N(d,s)=\frac{H_d}{H_{ds}}\Ex[\N(g,\Z_{s})]\quad\text{and}\quad R^{GCP}_\ell(d,s)=\frac{d+1}{ds+1}\Ex[\ell_t(g,\Z_{s})],
\end{equation*}where $g$ and $t$ are selected uniformly at random from the sets $\Ap(\Z_s)$ and $\Z_s$ respectively. 

For a fixed $s$, when $d$ (and consequently $q$) tend to infinity, we have that the asymptotic behavior of the ratios is given by the following. We use the classical approximation for the $d$th harmonic number given by $H_n=\log n+\gamma+\O(n^{-1})$, where $\gamma$ is the Euler-Mascheroni constant.

\begin{theorem}\label{smalls}
 Let $s$ be a fixed positive integer. If $d$ tends to infinity through an unbounded sequence for which $q=ds+1$ is a prime power, then the asymptotic behavior of  $R^{GCP}_\N(d,s)$ and $R^{GCP}_\ell(d,s)$ are given by 
    \begin{equation*}
    \lim_{d\to\infty}R^{GCP}_\N(d,s)=\Ex[\N(g,\Z_{s})]\quad\text{and}\quad \lim_{d\to\infty}R^{GCP}_\ell(d,s)=\frac{1}{s}\Ex[\ell_t(g,\Z_{s})].
\end{equation*}
\end{theorem}
\begin{proof}
 From Theorem \ref{expecGCP} we have that \begin{align*}
    \lim_{d\to\infty}R^{GCP}_\N(d,s)=\lim_{d\to\infty}\frac{\log d+\gamma+\O(d^{-1})}{\log d+\log s+\gamma+\O(d^{-1})}\Ex[\N(g,\Z_{s})]=\Ex[\N(g,\Z_{s})] 
\end{align*} and 
\begin{equation*}
    \lim_{d\to\infty}R^{GCP}_\ell(d,s)=\lim_{d\to\infty}\frac{d+1}{ds+1}\Ex[\ell_t(g,\Z_{s})]=\frac{1}{s}\Ex[\ell_t(g,\Z_{s})].
\end{equation*}
\end{proof}

Since $\Ex[\N(g,\Z_{s})]>1$ if $s>1$, we have that generalized cyclotomic permutations have asymptotically more cycles, on average, than random permutations. Increasing the index does not eliminate the arithmetic influence when $s$ is fixed and even as $d$ tends to infinity, the ratios of the expectations do not converge to $1$ in general.

We now provide the asymptotic behavior when $s=p^e$ using the results from Section \ref{SAff}. 

\begin{theorem}
    Let $d$ be a fixed positive integer and $p$ a prime. Along any unbounded sequence of positive integers $e$ for which $q=dp^e+1$ is a prime power, the asymptotic behavior of  $R^{\Gc}_\N(d,p^e)$ and $R^{\Gc}_\ell(d,p^e)$ are given by 
    \begin{equation*}
    R^{GCP}_\N(d,2^e)\sim\frac{H_d}{8\log 2}e\quad\text{and}\quad R^{GCP}_\ell(d,2^e)\sim\frac{4(d+1)}{9d}
\end{equation*} when $p=2$ and 
\begin{equation*}
    R^{GCP}_\N(d,p^e)\sim\frac{H_d(p-1)F}{2p\log p}e\quad\text{and}\quad R^{GCP}_\ell(d,p^e)\sim\frac{(d+1)}{d}\frac{p(p+1+G)}{(p-1)(p+1)^2}
\end{equation*} when $p$ is odd, where  $F=\displaystyle\sum_{\substack{r\mid p-1\\r\neq1}}\frac{\varphi(r)}{r}$, $G=\displaystyle\sum_{\substack{r\mid p-1\\r\neq1}}r\varphi(r)$ and $H_d$ is the $d$th harmonic number.
\end{theorem}
\begin{proof}
    We have that $\Ex[\N(\sigma)]=e\log p+\log d+\gamma+\O(p^{-e})$. Thus, the case $p=2$ follows directly from Theorems \ref{expeccycnumber2}, \ref{expeccyclength2} and \ref{expecGCP}. The case $p$ odd follows from Theorems \ref{expeccycnumberp}, \ref{expeccyclengthp} and \ref{expecGCP}.
\end{proof}

We observe that the expected number of cycles of GCPs of index $\frac{q-1}{p^e}$ grows quadratically in $e$ whereas in random permutations it only grows linearly. The expected cycle length under a GCP of a random element is asymptotically proportional to a random permutation, with a factor depending on $d$ and the divisors of $p-1$. 




\section{Conclusion and future research}\label{Con}
In this paper, we investigated two statistical properties of the functional graphs of affine permutations over $\Z_{p^e}$ and generalized cyclotomic permutations over finite fields: the expected number of cycles and the expected length of the cycle containing a randomly selected element. 
We established exact expressions for these expectations for generalized cyclotomic permutations of index $d$ in terms of the expectations for the underlying affine maps. In particular, we showed that the expected number of cycles is the product of the corresponding expectation for affine permutations and the $d$th harmonic number, while the expected length of the cycle containing a random element is multiplied by $\frac{(d+1)}{2}$. 

We analyzed the ratio between these statistical properties and those of uniformly random permutations on sets of the same cardinality. For affine permutations over $\Z_{p^e}$, the asymptotic results showed that the ratio of the expected number of cycles grows linearly with $e$ but tends to zero when $p$ grows through a sequence of safe primes. On the other hand, the ratio of the expected lengths of the cycles containing uniformly selected elements tends to a constant that depends on $p$, as $e$ grows and approaches $\frac{3}{2}$ when $p$ grows through a sequence of safe primes.

In generalized cyclotomic permutations, the asymptotic results show that for $s$ fixed, as $d$ tends to infinity the ratios of the expectations approach constants determined by the corresponding expectations for affine permutations over $\Z_s$. On the other hand, when $s=p^e$ with $p$ and $d$ fixed, the expected number of cycles grows quadratically with $e$, whereas the corresponding expectation for uniformly random permutations grows linearly. Their expected cycle lengths for a random element remain asymptotically proportional, with a factor depending on $p$ and $d$.

An interesting direction for future research is the extension of the explicit expectation formulas to affine permutations over $\mathbb Z_n$ for arbitrary positive integers $n$. Although their functional graphs can be determined from the prime decomposition of $n$, the interaction between cycle lengths from different components makes the calculation of these expectations more involved. Another possible direction is the study of higher moments and the distribution of cycle lengths, which could provide a more detailed understanding of how the functional graphs of these structured permutation families differ from those of uniformly random permutations.

\section*{Acknowledgments}
The author was funded by the Austrian Science Fund [10.55776/PAT4719224].

\bibliographystyle{plain}
\bibliography{Bibliography}
\end{document}